\documentstyle[epsf]{article}

\font\myfont=msbm10 scaled \magstep1
\def\CC{\hbox{\myfont\char'103}}

\def\ss{s_n(R_n z, E_{\frac{1}{\rho}})}
\def\tt{t_{n+1}(R_n z, E_{\frac{1}{\rho}})}
\def\GG{\Gamma\left(1+\frac{n}{\rho}\right)}
\def\G{\Gamma\left(1-\frac{1}{\rho}\right)}

\title{Linear combinations of sections and tails of
Mittag--Leffler functions and their zeros}
\author{N.A. Zheltukhina}

\begin{document}

\maketitle

\noindent{\bf Abstract}. The zero distribution of sections of
Mittag--Leffler functions of order $\rho>1$ was studied in 1983 by
A.~Edrei, E.B.~Saff and R.S.~Varga. In the present paper, we study
the zero distribution of linear combinations of sections and tails
of Mittag--Leffler functions of order $\rho >1$.

\bigskip

\bigskip

\renewcommand{\baselinestretch}{1.4}

\begin{large}

\centerline{\bf 1. Introduction and statement of results}

\bigskip

Let
$$
f(z)=\sum_{k=0}^\infty a_k z^k, \quad a_0>0,
\eqno(1.1)
$$
be a transcendental  entire function. Denote by
$$
s_n(z,f)=\sum_{k=0}^n a_k z^k
\eqno(1.2)
$$
and
$$
t_n(z,f)=\sum_{k=n}^\infty a_k z^k
\eqno(1.3)
$$
its sections and tails respectively.
For $r>0$, set
$$
\mu(r):=\max_{k} |a_k|r^k, \qquad
\nu(r):=\max\{k:|a_k|r^k=\mu(r)\}.
\eqno(1.4)
$$
The function $\nu(r)$, called the central index,
is an increasing step function of $r$ taking
integer values ( see \cite[p.3]{[28]}). Denote
by $R_1,R_2,R_3,\ldots$, the discontinuity points of
$\nu(r)$. Note that
$\nu(r)$ is everywhere continuous from the right,
and $\lim\limits_{n\to\infty } R_n=\infty$ (see \cite[pp.5--6]{[28]}).

Let ${\cal{M}}_n(\lambda,f)$, $\lambda\in\CC$, be the set of all roots of the equation
$$
(1-\lambda)s_n(R_nz,f)-\lambda t_{n+1}(R_nz,f)=0.
\eqno(1.5)
$$
In particular, ${\cal{M}}_n(0,f)$
coincides with the set of zeros of the sections (1.2)
and ${\cal{M}}_{n-1}(1,f)$ coincides with the set of zeros of the tails (1.3).
Define ${\cal{M}}(\lambda,f)$ to be the set of all accumulation points of
$\cup_{n=1}^\infty {\cal{M}}_n(\lambda,f)$.

In 1924, G. Szeg\"o \cite{[17]} proved a remarkable theorem related to the
behavior of the roots of the equation (1.5) for $f(z)=e^z$. For function
$f(z)=e^z$ we have $R_n=n$ (see \cite[Problem 27, p.6]{[28]}). Consider
the main result of \cite{[17]}.
The curve
$$
S:=\{z:\quad |ze^{1-z}|=1\}
$$
is known as the Szeg\"o curve.


\bigskip

\noindent {\bf Theorem S} (\cite{[17]}). {\it The following equality holds}
$$
{\cal{M}}(\lambda, e^z)=\left\{\begin{array}{ll}
S\cap\{z: |z|\leq 1\},& for \quad \lambda=0;\\
S\cap\{z: |z|\geq 1\}
,& for \quad \lambda=1;\\
S,& for \quad \lambda\ne 0,1.
\end{array}\right.
\eqno(S)
$$

In \cite{[4]}, A. Edrei, E.B. Saff and R.S. Varga  investigated the
distribution of the zeros of the sections (1.2) of
Mittag--Leffler functions of all
orders $\rho>1$:
$$
E_{\frac{1}{\rho}}(z)=\sum_{j=0}^\infty \frac{z^j}{\Gamma(1+
\frac{j}{\rho})},\qquad 1<\rho<\infty.
\eqno(1.6)
$$
For function $E_{\frac{1}{\rho}}(z)$, by \cite[Problem 117, p.26]{[29]},
we have
$$
R_n=\displaystyle{
\frac{\Gamma\left(1+\frac{n}{\rho}\right)}
{\Gamma\left(1+\frac{n-1}{\rho}\right)}}.
$$
It follows from Stirling's formula for $\Gamma (x)$, $x>0$, that
$$
R_n= \left(\frac{n}{\rho}\right)^{
\frac{1}{\rho}}\left(1+\frac{\rho-1}{2\rho
n}+O\left(\frac{1}{n^2}\right) \right), \quad n\to\infty.
\eqno(1.7)
$$
Consider the main result of \cite{[4]}.
By $S(\rho)$ denote a  curve which is
represented, in polar coordinates, as the set of all points
$$
\psi(\phi)=\sigma(\phi)e^{i\phi}, \qquad
-\frac{\pi}{2\rho}\leq\phi<2\pi-\frac{\pi}{2\rho},
$$
where\\
(i) for
  $-\displaystyle{\frac{\pi}{2\rho}}\leq \phi\leq
\displaystyle{\frac{\pi}{2\rho}}$,
$\sigma(\phi)$ satisfies
$$
\{\sigma(\phi)\}^\rho \cos(\phi \rho) - 1-\rho\log \sigma(\phi)=0,
\eqno(1.8)
$$
(ii) for
  $\displaystyle{\frac{\pi}{2\rho}}< \phi<
2\pi -\displaystyle{\frac{\pi}{2\rho}}$,
we take
$$
\sigma(\phi)=e^{-\frac{1}{\rho}}.
$$
The curve $S(\rho)$ is called the generalized Szeg\"o curve.
We will prove in Lemma 1 that
the curve $S(\rho)\cap \{z: |z|\geq 1\}$
 has asymptotes
$\arg z=\pm\displaystyle{\frac{\pi}{2\rho}}$
( meanwhile, the original Szeg{\"o} curve $S$ doesn't have asymptotes ).

The arguments of the zeros of $E_{\frac{1}{\rho}}(z)$ for
$1<\rho<\infty$ tend to
$\displaystyle{\pm \frac{\pi}{2\rho}}$ as $|z|\to\infty$.
Hence, as a consequence of the uniform convergence
$$
s_n(z, E_{\frac{1}{\rho}})\to
 E_{\frac{1}{\rho}}(z), \quad n\to\infty,
$$
on compact subsets, there are zeros of $s_n(R_nz, E_{\frac{1}{\rho}})$ whose
arguments are close to
$\displaystyle{\pm \frac{\pi}{2\rho}}$.
Denote
$${\cal{M}}_0(\lambda,E_{\frac{1}{\rho}})={\cal{M}}
(\lambda, E_{\frac{1}{\rho}})
\backslash \left\{z:\arg z=\pm \frac{\pi}{2\rho}\right\}.
$$
The main result of \cite{[4]}
implies that ${\cal{M}}_0(0,E_{\frac{1}{\rho}})=S(\rho)\cap\{z:|z|\leq 1\}$.
The question arises whether the complete analogue of the Szeg{\"o}
theorem holds for $E_{\frac{1}{\rho}}(z)$, that is whether the equality
holds
$$
{\cal{M}}_0(\lambda, E_{\frac{1}{\rho}})=\left\{\begin{array}{llll}
S(\rho)\cap\{z: |z|\leq 1\},& for \quad \lambda=0;&
\qquad \qquad\qquad \qquad \qquad&(S')\\
S(\rho)\cap\{z: |z|\geq 1\}
,& for \quad \lambda=1;&\qquad&(S'')\\
S(\rho),& for \quad \lambda\ne 0,1.&&(S''')
\end{array}\right.
$$
In this paper, we are going to prove $(S'')$ and $(S''')$, and
extensions of main results of \cite{[4]} for the linear
combination $(1-\lambda) s_n(R_n z,E_{\frac{1}{\rho}})- \lambda
t_{n+1}(R_n z,E_{\frac{1}{\rho}})$, $\rho >1$. The main result of
this paper is the following theorem dealing with the asymptotic
expressions of this linear combination in different domains of
$\CC$.

\bigskip

\noindent{\bf Theorem 1}. {\it Let
$$
I_n(R_n z;\lambda):=(1-\lambda)\ss -\lambda\tt,
$$
$\delta_2$, $\delta_3$ be given positive constants, and $\rho >1$.
Then, as $n\to\infty$,
$$
\frac{I_n(R_n z;\lambda)\GG}{R_n^n z^n}=
-\lambda \rho\frac{e^{R_n^\rho z^\rho}\GG}{R_n^n z^n}(1+o(1))
-\frac{z}{1-z}(1+o(1)),
$$
$$
\mbox{if}\quad z\in\{z=r e^{i\phi}: r\geq 1,|\phi|\leq \frac{\pi}{2\rho},
 |z-1|\geq \delta_2\},
$$
$$
\frac{I_n(R_n z;\lambda)\GG}{R_n^n z^n}=
(1-\lambda) \rho\frac{e^{R_n^\rho z^\rho}\GG}{R_n^n z^n}(1+o(1))-
\frac{z}{1-z}(1+o(1)),
$$
$$
 \mbox{if}\quad z\in\{z=r e^{i\phi}: 0< r\leq 1,
|\phi|\leq \frac{\pi}{2\rho}, |z-1|\geq \delta_2\},
$$
$$
\frac{I_n(R_n z;\lambda)\GG}{R_n^n z^n}=
\frac{(\lambda-1)}{\G} \frac{\GG}{ R_n^{n+1} z^{n+1}}(1+o(1))
-\frac{z}{1-z}(1+o(1)),
$$
$$
\mbox{if}\quad z\in\{z=r e^{i\phi}: r> 0,|\phi|\geq
\frac{\pi}{2\rho}+
\delta_3,\}.
$$
In all expressions above, the symbols $o(1)$ have uniformity property
with respect to $z$. }

We remark that, in the special case $\lambda =1$, the asymptotic
expression for $I_n(R_n z; 1)$ in the unit disc, one can find in
\cite[Lemma 9.2]{[4]}. It coincides with one we obtain in Theorem
1.
\bigskip

For given  $\delta_2>0$, $\delta_3>0$ and $h>0$, let us
introduce the following regions.
$$
\begin{array}{l}
 \Omega_1= \{z=re^{i\phi}: 0< r\leq 1, |z-1|\geq \delta_2,
|\phi|\leq
\frac{\pi}{2\rho}-\delta_3,  r^\rho\cos(\rho\phi)-1-\rho\log r\geq 0\},\\
\Omega_2= \{z=re^{i\phi}: |\phi|\leq
\frac{\pi}{2\rho}-\delta_3,  r^\rho\cos(\rho\phi)-1-\rho\log r\leq -h\},\\
\Omega_3=\{z=r e^{i\phi}: e^{-1/\rho}+h\leq r,  |\phi|\geq
\frac{\pi}{2\rho}+\delta_3\},\\
\Omega_4=\{z=re^{i\phi}: 0<r\leq e^{-1/\rho}-h,  |\phi|\geq
\frac{\pi}{2\rho}+\delta_3\},\\
\Omega_5=\{z=re^{i\phi}:r\geq 1, |\phi|\leq
\frac{\pi}{2\rho}-\delta_3,  r^\rho\cos(\rho\phi)-1-\rho\log r\geq h\}.\\
\end{array}
$$


\noindent The next theorem deals with the zero free regions of $I_n(R_n z;
\lambda )$.

\bigskip

\noindent{\bf Theorem 2}. {\it Let $I_n(R_n z;\lambda)$ be as in Theorem 1,
  $\delta_2$, $\delta_3$ and $h$ be
given  positive constants.
 Then, for all sufficiently large $n$,
 $I_n(R_nz;\lambda)$ has no zeros in $\cup_{i=1}^5\Omega_i$.}

\bigskip

\noindent Theorem 2 is an extension  of Theorem 5 from \cite{[4]}.
The next two theorems give  information on the zero distribution
of $I_n(R_n z;\lambda)$ in the neighborhood of points on the
generalized Szeg\"o curve $S(\rho)$. To characterize this
distribution in the neighborhood of the point $z=1$, denote by
$$
\mbox{erfc}(\zeta)=1-\frac{2}{\pi}\int_{0}^\zeta e^{-v^2}dv.
$$
 the complementary error function.

\bigskip

\noindent{\bf Theorem 3}. {\it With $I_n(R_n z;\lambda)$
defined in Theorem 1, we have}
$$
\left(1+\left(\frac{2}{\rho n}\right)^{1/2}\zeta\right)^{-n}
\left\{E_{\frac{1}{\rho}}(R_n)\right\}^{-1}
I_n\left(R_n\left(1+\left(\frac{2}{\rho n}\right)^{1/2}\zeta\right); \lambda
    \right)
\to e^{\zeta^2}\left\{\frac{\mbox{erfc}(\zeta)}{2}-\lambda\right\}, \quad
$$
{\it as $n\to\infty$,
uniformly on every compact set of the $\zeta$--plane.}

\bigskip

\noindent Theorem 3  is an extension  of Theorem 1 from \cite{[4]}.

\bigskip

\noindent{\bf Theorem 4}.\\
I. {\it Let $\xi=\xi(\phi)$, $0<\phi<\frac{\pi}{2
\rho}$, be a fixed point on the generalized Szeg{\"o} curve $S(\rho)$. Let
$$
\tau=|\zeta|^\lambda \sin(\phi \rho) -\rho\phi,
$$
and let the sequence $\{\tau_n\}_{n}$ be defined by the condition
$$
\tau_n\equiv
\frac{\tau}{\rho}n(\bmod\, 2\pi), \qquad -\pi<\tau_n\leq\pi.
$$
Then, as $n\to\infty$,
$$
I_n\left(R_n \xi \left(1+\frac{\log n}{2(1-\xi^\rho)n}-
\frac{\zeta- i\tau_n}{(1-\xi^\rho)n}\right);\lambda\right)
\displaystyle{ \frac{\GG}{R_n^n \xi^n\left(1+\frac{\log
n}{2(1-\xi^\rho)n}- \frac{\zeta- i\tau_n}{(1-\xi^\rho)n}\right)^n}
} \eqno(1.9)
$$
$$
\to
\left\{\begin{array}{ll}
\displaystyle{
(1-\lambda)(2\pi\rho)^{\frac{1}{2}}e^{\frac{\rho+1}{2\rho}(\xi^\rho-1)} e^{
\zeta}-\frac{\xi}{1-\xi}}, & \mbox{if}\quad  |\xi|< 1,\\
-\displaystyle{
\lambda(2\pi\rho)^{\frac{1}{2}}e^{\frac{\rho+1}{2\rho}(\xi^\rho-1)} e^{
\zeta}-\frac{\xi}{1-\xi}}, & \mbox{if}\quad  |\xi|>1,
\end{array}\right.
$$
 uniformly on every compact set of the $\zeta$ -- plane.

\noindent II. Let $\xi=e^{-\frac{1}{\rho}}e^{i\phi}$,
$\frac{\pi}{2\rho} < \phi\leq \pi$, be a fixed point on the
circular portion of $S(\rho)$, and let the sequence $\tau'_n$ be
defined by the condition
$$
\tau^{'}_n\equiv(n+1)\phi (mod\quad 2\pi), \quad -\pi<\tau'_n\leq\pi.
$$
Then
$$
I_n\left(R_n \xi \left(1+\left(\frac{1}{2}-\frac{1}{\rho}\right)
\frac{\log n}{n}- \frac{\zeta-
i\tau^{'}_n}{n+1}\right);\lambda\right) \displaystyle{
\frac{\GG}{R_n^n \xi^n\left(1+(\frac{1}{2}-\frac{1}{\rho})
\frac{\log n}{n}- \frac{\zeta- i\tau^{'}_n}{n+1}\right)^n}}
\eqno(1.10)
$$
$$
\to \frac{(\lambda-1)(2\pi
e^{\frac{1-\rho}{\rho}})^{\frac{1}{2}}}{ \rho^{\frac{1}{2}-
\frac{1}{\rho}}\G} e^{- \zeta}-\frac{\xi}{1-\xi}
$$
uniformly on every compact set of the $\zeta$ -- plane.}

\bigskip

\noindent Arguments of $I_n$ in (1.9) and (1.10) of Theorem 4 were
taken from analogous Theorems 2 and 3 from
\cite{[4]}. Theorem 4 is an easy corollary of  Theorem 1.

\bigskip

\centerline{\bf 2.   Preliminary results }

\bigskip

\noindent{\bf Lemma 1.} {\it The curve $S(\rho)\cap \{z: |z|\geq 1\}$
 has asymptotes $l_j=
\left\{\arg z=(-1)^{j-1}\displaystyle{\frac{\pi}{2\rho}}\right\}$, $j=1,2$.}

\bigskip

\noindent{\it Proof }.
 It suffices to consider only the upper half $B(\rho)$ of the curve.
It is given by the equation
$$
\cos \rho \phi =r^{-\rho}(1+\rho \log r),
$$
which implies  that conditions
$r\to\infty$ and $\phi\to \displaystyle{\frac{\pi}{2\rho}}$ are equivalent.
Take a point
$A=(r\cos \phi,r\sin\phi)$ on $B(\rho)$
and a point $C=\left(r\cos\phi,
r\tan\displaystyle{\frac{\pi}{2\rho}}\cos\phi\right)$ on $l_1$.
The distance $d(A,l_1)$  from point $A$ to $l_1$ does not exceed the distance
between points $A$ and $C$, i.e.
$$
d(A, l_1)\leq r\tan \frac{\pi}{2\rho} \cos \phi -r\sin \phi=
\leq \frac{r\sin \rho\left(\displaystyle{\frac{\pi}{2\rho}}-
\phi\right)}{\cos\displaystyle{\frac{\pi}{2\rho}}}=
\frac{1+\rho \log r}{r^{\rho -1}\cos
\displaystyle{\frac{\pi}{2\rho}}}\to 0 \quad {\mbox as} \quad r\to\infty.
$$
$\Box$

\bigskip

With $I_n(R_n z;\lambda)$ defined in the statement of
Theorem 1, rewrite (1.5) as
$$
I_n(R_nz;\lambda)=(1-\lambda)E_{\frac{1}{\rho}}(R_n z)-\tt=0.
\eqno(2.1)
$$
In the proof of Theorem 1, by (2.1), asymptotic formulas for
$E_{\frac{1}{\rho}}(R_n z)$ and $\tt$ will play a major role. Let
us present those we will use.

The functions $E_{\frac{1}{\rho}}(z)$ have the following
well--known asymptotic relations
(see \cite[p.114]{[24]}):
$$
E_{\frac{1}{\rho}}(z)=\left\{
\begin{array}{lll}
\rho e^{z^{\rho}}-\displaystyle{\frac{1}{z\G}}+O
\displaystyle{\left(\frac{1}{|z|^2}\right)},
\quad& |\arg z|\leq \displaystyle{\frac{\pi}{2\rho}},& |z|\to\infty,\\
-\displaystyle{\frac{1}{z\G}}+O
\displaystyle{\left(\frac{1}{|z|^2}\right)}, \quad &
\displaystyle{\frac{\pi}{2\rho}}\leq |\arg z|\leq \pi, &|z|\to\infty .
\end{array}\right.
\eqno(2.2)
$$

To write an asymptotic expression for $\tt$,
we will use so-called  Mittag--Leffler type functions
$$
E_{\frac{1}{\rho}}(z,\mu)=\sum_{k=0}^\infty  \displaystyle{
\frac{z^k}{\Gamma\left(\mu+
\frac{k}{\rho}\right) }}, \qquad \rho>0, \quad \mu\in \CC,
$$
introduced by M.M. Djrbashian (see \cite[p.117]{[27]}).
Denote by $L(\alpha, H)$ ($H>0$, $0<\alpha\leq\pi$) a contour
following nondecreasing direction of $\arg \zeta$
and consisting of  two rays $\arg\zeta=\pm\alpha$, $|\zeta|\geq H$
and an arc $-\alpha\leq \arg \zeta\leq\alpha$ of a circle $|\zeta|=H$.
By $G^-(H,\alpha)$ and $G^+(H, \alpha)$ we denote  two regions lying
respectively from the left and right sides of $L(\alpha, H)$.
For $\rho>1$ and $\pi/(2\rho)<\nu\leq \pi/\rho$,
the following  representations hold
(\cite[p.127]{[27]})
$$
E_{\frac{1}{\rho}} (z,\mu)=
\left\{\begin{array}{ll}
\rho z^{\rho(1-\mu)} e^{z^\rho}
+ \displaystyle{\frac{\rho}{2\pi i }} \int\limits_{L(\nu,H)}
\displaystyle{\frac{e^{\zeta^\rho} \zeta^{\rho(1-\mu)}}{\zeta -z}} d\zeta,
&z\in G^+(H,\nu),\\
\displaystyle{\frac{\rho}{2\pi i }} \int\limits_{L(\nu,H)}
\displaystyle{\frac{e^{\zeta^\rho} \zeta^{\rho(1-\mu)}}{\zeta -z}} d\zeta,&
z\in G^-(H,\nu).
\end{array}\right.
$$
Since
$$
\tt=\sum_{k=n+1}^{\infty}\frac{(R_nz)^k}{\Gamma (1+\frac{k}{\rho})}=
(R_nz)^{n+1}E_{\frac{1}{\rho}}\left(R_nz,1+\frac{n+1}{\rho}\right),
$$
we have
$$
\tt=
\left\{
\begin{array}{l}
\displaystyle{\rho e^{R_n^\rho z^\rho}
 +
\frac{\rho(R_n z)^{n+1}}{2\pi i}\int\limits_{L(\nu,H)}
\frac{e^{\zeta^\rho} \zeta^{-(n+1)}}{\zeta-R_nz} d\zeta},
\quad R_nz\in G^+(H,\nu),\\
\frac{\rho(R_n z)^{n+1}}{2\pi i}\int\limits_{L(\nu,H)}
\displaystyle{
\frac{e^{\zeta^\rho} \zeta^{-(n+1)}}{\zeta -R_nz} d\zeta},\quad
R_n z\in G^-(H,\nu),
\end{array}\right.
$$
where $\nu$ and $H>0$ are any constants but
$\displaystyle{\frac{\pi}{2\rho}}<\nu\leq \frac{\pi}{\rho}$.
Then it follows from (2.1) and (2.2) that, as $n\to\infty$,
$$
\begin{array}{llll}
I_n(R_nz;\lambda)=\quad \displaystyle{-\lambda\rho e^{R_n^\rho z^\rho}(1+o(1))}
&-&
\displaystyle{
\frac{\rho(R_n z)^{n+1}}{2\pi i}\int\limits_{L\left(\frac{\pi}{2\rho}+\frac{
\delta_3}{2},R_n\right)}
\frac{e^{\zeta^\rho} \zeta^{-(n+1)}}{\zeta-R_nz} d\zeta},
&\qquad (2.3)\\
 \mbox{if} \quad |z|>1 \quad \mbox{and} \quad |\arg z|\leq
\frac{\pi}{2\rho},& &\\
I_n(R_nz;\lambda)=\displaystyle{
\frac{(\lambda -1)}{R_n z \Gamma\left(1-\frac{1}{\rho}\right)}
(1+o(1))}&-&
\displaystyle{\frac{\rho(R_n z)^{n+1}}{2\pi i}
\int\limits_{L\left(\frac{\pi}{2\rho}+\frac{
\delta_3}{2},R_n\right)}
\frac{e^{\zeta^\rho} \zeta^{-(n+1)}}{\zeta -R_nz} d\zeta},&\qquad (2.4)\\
\mbox{if}\quad  |z|>0 \quad \mbox{and} \quad |\arg z|>
\frac{\pi}{2\rho}+\frac{\delta_3}{2\rho},& &\\
I_n(R_nz;\lambda)=\displaystyle{(1-\lambda)\rho e^{R_n^\rho z^\rho}(1+o(1))}&
-&\displaystyle{
\frac{\rho(R_n z)^{n+1}}{2\pi i}\int\limits_{L\left(\frac{\pi}{2\rho}+\frac{
\delta_3}{2},R_n\right)}
\frac{e^{\zeta^\rho} \zeta^{-(n+1)}}{\zeta-R_nz} d\zeta},
&\qquad(2.5)\\
\mbox{if}\quad 0<|z|<1 \quad \mbox{and} \quad
|\arg z|\leq \frac{\pi}{2\rho}.& &\\
\end{array}
$$
Here and in the text below estimates denoted with symbols $o$ and $O$ are
uniform with respect to $z$.

\bigskip

\centerline{\bf 3. Proof of Theorem 1 }

\bigskip

Set
$$
K_n(z):=\displaystyle{\int\limits_{L\left(\frac{\pi}{2\rho}+\frac{
\delta_3}{2},R_n\right)} \frac{e^{\zeta^\rho}
\zeta^{-n}}{\zeta-R_nz} d\zeta}.
\eqno(3.1)
$$
To prove Theorem 1, by (2.3) -- (2.5) and (3.1),
  it is enough
to find  asymptotic expression for the integral $K_n(z)$. We will do this in
the following three steps:\\
1) change the contour of integration of $K_n(z)$;\\
2) show that the main contribution to $K_n(z)$ comes from the neighborhood of the point $\zeta=R_n$\\
3) find asymptotic expression for $K_n(z)$ by using Laplace's Method
for contour integrals.

To introduce new contours of integration we define the following curves.
Denote by $T(\rho)$ a simple curve which is represented in polar coordinates as the set of all points
$$
r^{\rho}=\frac{\rho \phi}{\sin \rho\phi}, \qquad -\frac{\pi}{\rho}<\phi<
\frac{\pi}{\rho}.
$$
Define
$$
S(\rho,h):=\left\{z=re^{i\phi}:  r^\rho\cos(\rho\phi)-
\rho\log r-1=-\frac{h}{2}\right\}.
$$

\noindent By the same reason as in Lemma 1, $S(\rho,h)$ has two
asymptotes $\arg z=\pm\displaystyle{\frac{\pi}{2\rho}}$. The curves $S(\rho,
h)$ and $T(\rho)$ have two common points,
$z_1=de^{i\gamma}$ and $z_2=de^{-i\gamma}$, say.
 Note that constants $d$ and
$\gamma$ depend only on $h$, and as $h\to 0$ we see that
$\gamma\sim\displaystyle{
\frac{\sqrt{h}}{\rho}}$ and $d^\rho\sim\displaystyle{
\frac{\sqrt{h}}{\sin\sqrt{h}}}$.
Denote by $TS_1(\rho,h)$ the curve which consists of\\
$(i)$ $T(\rho)\cap\{z:|z|\leq d\}$, where $d=d(h)$
is defined above;\\
$(ii)$ $  S(\rho,h)\cap\left\{z:|z|\geq d,
|\arg z|\leq \frac{\pi}{2\rho}-\delta_3\right\}$;\\
$(iii)$ $\{z:|z|=D_1, \frac{\pi}{2\rho}-\delta_3\leq |\arg z|\leq
\frac{\pi}{2\rho}+\frac{\delta_3}{2}\}$,
 where $D_1$ is such
that\\
$D_1^{\rho}\sin(\rho \delta_3)-\rho\log D_1-1=-\displaystyle{\frac{h}{2}}$;\\
$(iv)$ $\{z: \arg z=\pm\left(\frac{\pi}{2\rho}+
\frac{\delta_3}{2}                \right)
                             ,\quad  |z|\geq D_1\}$,
where $D_1$ was defined in part (iii).


\noindent Denote by $TS_2(\rho,h)$ the curve which consists of \\
$(i)$ defined above;\\
$(ii')$ $  S(\rho,h)\cap\left\{z:|z|\leq d, |\arg z|\leq
\frac{\pi}{2\rho}+\frac{\delta_3}{2}\right\}$, where $d$ is the
same as
in (i);\\
$(iii')$ $\{z: \arg z=\pm\left(\frac{\pi}{2\rho}+
\frac{\delta_3}{2}\right), |z|\geq D_2\}$, where $D_2$ is such
that\\
 $D_2^\rho\sin\left(\rho\frac{\delta_3}{2}\right)+\rho\log
D_2 +1=\frac{h}{2}$.

\bigskip

\noindent{\bf Lemma 2}.  {\it Let
$$
z\in \Omega :=\CC\backslash \left\{\quad\left\{z:\frac{\pi}{2\rho}\leq
|\arg z|\leq \frac{\pi}{2\rho}+\delta_3\right\}\cup\{z:|z-1|\leq \delta_2\}
\quad\right\}
$$
Then, for sufficiently small $h>0$,
$$
K_n(z)=\frac{
\displaystyle{e^{\frac{n}{\rho}}
           }}{R_n^{n+1}}\int\limits_{(i)}
\frac{ \displaystyle{e^{\frac{n}{\rho}(t^\rho-\rho\log t-1)}
e^{\frac{\rho-1}{2\rho}t^{\rho}}}
}{(t-z)t}\left(1+O\left(\frac{1}{n}\right)t^{\rho}\right)dt +
O\left(\frac{e^{\frac{n}{\rho}(1-\frac{h}{2})}}{
R_n^{n+1}}\right),
$$
where $K_n(z)$ is defined by (3.1).}

\bigskip

\noindent {\it Proof}. Let $z\in G^{-}(1, \frac{\pi}{2\rho}+
\frac{\delta_3}{2})\cap\Omega$.   For sufficiently small $h>0$,
using (1.7), we have
$$
K_n(z)=\frac{1}{R_n^{n+1}}\int\limits_{L\left(1,
\frac{\pi}{2\rho}+ \frac{\delta_3}{2}\right)
}\frac{e^{R_n^{\rho}t^{\rho}} t^{-(n+1)}}{t-z}dt
=\frac{1}{R_n^{n+1}}\int\limits_{TS_1(\rho,h)}
\frac{e^{R_n^{\rho}t^{\rho}} t^{-n}}{(t-z)t}dt
$$
$$
=\frac{1}{R_n^{n+1}}\int\limits_{TS_1(\rho,h)}
\frac{\displaystyle{
e^{\frac{n}{\rho}t^{\rho}} t^{-n} e^{\frac{\rho-1}{2\rho}t^{\rho}}}
}{(t-z)t}\left(1+O\left(\frac{1}{n}\right)t^{\rho}\right)dt.
$$
Since there exists such constant $A_1>0$ that $|t-z|>A_1$
for all $z\in G^{-}(1, \frac{\pi}{2\rho}+
\frac{\delta_3}{2})\cap \Omega$  while $t\in TS_1(\rho,h)$,
for any number  $l$ we have
$$
\left|\int\limits_{(ii)}
\frac{\displaystyle{
e^{\frac{n}{\rho}t^{\rho}} t^{l-n} e^{\frac{\rho-1}{2\rho}t^{\rho}}}
}{t-z}dt\right|\leq B_1\int\limits_{(ii)}
e^{\frac{n}{\rho}|t|^\rho\cos(\rho\phi)}|t|^{-n}d|t|=
O\left(e^{\frac{n}{\rho}(1-\frac{h}{2})}\right);
$$
$$
\left|\int\limits_{(iii)}
\frac{\displaystyle{
e^{\frac{n}{\rho}t^{\rho}} t^{l-n} e^{\frac{\rho-1}{2\rho}t^{\rho}}}
}{t-z}dt\right|\leq B_2
 \int\limits_{(iii)}
\displaystyle{
e^{\frac{n}{\rho}D_1^{\rho}\sin(\rho\delta_3)}} D_1^{-n} d|t|
=O\left(e^{\frac{n}{\rho}(1-\frac{h}{2})}\right);
$$
$$
 \left|\int\limits_{(iv)}
\frac{\displaystyle{
e^{\frac{n}{\rho}t^{\rho}} t^{l-n} e^{\frac{\rho-1}{2\rho}t^{\rho}}}
}{t-z}dt\right|\leq B_3
 \int\limits_{(iv)}
\displaystyle{
e^{-\frac{n}{\rho}|t|^{\rho}\sin\left(\frac{\rho\delta_3}{2}\right)}} d|t|
=o(1);
$$
where $B_1$, $B_2$, $B_3$ are constants not depending on $z$.
The last three inequalities prove the Lemma  if
$z\in G^{-}(1, \frac{\pi}{2\rho}+\frac{\delta_3}{2})\cap\Omega$.
The same proof works for $z\in G^{+}
(1, \frac{\pi}{2\rho}+\frac{\delta_3}{2})\cap\Omega$
 if we change the
contour of integration in the integral $K_n(z)$ by $TS_2(h,\rho)$.
It follows that Lemma 2 holds  for $z\in \Omega$.
$\Box$

\bigskip

\noindent{\bf Lemma 3}. {\it Let $|z-1|\geq \delta_2$
and $l$
be a fixed real number.
Then, for sufficiently small $h>0$,}
$$
\int\limits_{(i)}
\frac{\displaystyle{
e^{\frac{n}{\rho}(t^{\rho}-\rho\log t-1)}  e^{\frac{\rho-1}{2\rho}t^{\rho}}}
t^l}{t-z}dt=
-\frac{i\sqrt{2\pi}e^{\frac{\rho-1}{2\rho}}(1+o(1))}{\rho
(1-z)\left(\frac{n}{\rho}
\right)^{1/2}}.
$$

\bigskip

\noindent{\it Proof}. Rewrite
$$
\int\limits_{(i)}
\frac{\displaystyle{
e^{\frac{n}{\rho}(t^{\rho}-\rho\log t-1)}  e^{\frac{\rho-1}{2\rho}t^{\rho}}}
t^l}{t-z}dt=
\int_{{\cal{D}}} e^{-\frac{n}{\rho}v}f(v)dv,
$$
where $v=-t^\rho +\rho\log t+1$, $f(v)dv= \displaystyle{\frac{
e^{\frac{\rho-1}{2\rho}t^{\rho}}t^l}{t-z}}dt$, or
equivalently, $ f(v)= \displaystyle{\frac{
e^{\frac{\rho-1}{2\rho}t^{\rho}}t^{l+1}}{\rho(t-z)(1-t^\rho)}}$,
and ${\cal{D}}={\cal{D}}_1\cup{\cal{D}}_2$, where ${\cal{D}}_1$
is the upper side of
the segment $[0;d_1]$ following the  direction of the decrease of  $v$ and
${\cal{D}}_2$ is the lower side of the segment $[0;d_1]$
following the  direction of increase of  $v$.
 Here $d_1=\displaystyle{ \frac{h}{2}=-d^\rho \cos (\rho\gamma)
+\rho\log d+1}$, where $d$ and $\gamma$  are the polar coordinates
of a point of intersection of $T(\rho)$ and $S(\rho,h)$.

Note  that the function $v=-t^\rho+\rho\log t+1$ maps the region on the
right side of $T(\rho)$ conformally onto the whole $v$--plane cut
 along the positive ray. In particular, the image of the curve
$(i)$ is the segment $[0;d_1]$ traced twice. The transformation
$w=\sqrt{v}$ maps $v$--plane  cut along the positive ray onto the
upper half plane. We have $w^2=-t^\rho+\rho\log t+1=-
\displaystyle{\frac{{\rho}^2}{2}(t-1)^2}\psi(t)$, where $\psi(t)$
is an analytic function in some neighborhood of $t=1$ and $\psi
(1)=1$. Then  $ w=\displaystyle{\frac{\rho}{
\sqrt{2}}i(t-1)}\psi_1(t)$, where $\psi_1(t)$ is an analytic
function  in some neighborhood of $t=1$ and $\psi_1(1)=1$. Since
$w$ is analytic in a neighborhood of $t=1$ and $w'(1)=
\displaystyle{\frac{\rho i}{\sqrt{2}}}\ne 0$, then its inverse
function $t(w)$ is analytic in a neighborhood of $w=0$, and hence
the following function
$$
g(w):=wf(w^2)=\frac{ i (t-1)
\psi_1(t)e^{\frac{\rho-1}{2\rho}t^\rho}t^{l+1}}{
\sqrt{2}(t-z)(1-t^\rho)}=-\frac{ i \psi_1(t)e^{
\frac{\rho-1}{2\rho}t^\rho}t^{l+1}}{
\sqrt{2}\rho(t-z)}(1+o(1)), \quad |t|\to 1,
$$
 is analytic in some neighborhood of $w=0$, say $
|w|\leq C$, where $C$ is  a constant not depending on $z$. Let $|w|<C/2$, then
$$
g(w)=\frac{1}{2\pi i} \int\limits_{|\zeta|=C}\frac{g(\zeta)}{\zeta-w}d\zeta=
\frac{1}{2\pi i}\int\limits_{|\zeta|=C}\frac{g(\zeta)}{\zeta}d\zeta+
\frac{w}{2\pi i}\int\limits_{|\zeta|=C}\frac{g(\zeta)}{\zeta(\zeta-w)}d\zeta
$$
$$
=g(0)+w\alpha(w)=-\frac{ i e^{
\frac{\rho-1}{2\rho}}}{
\sqrt{2}\rho(1-z)}+w\alpha(w),
$$
where $\alpha(w)$ is a function analytic in $|w|<C/2$, and
$$
|\alpha(w)|\leq \frac{2\pi C\max\limits_{|\zeta|=C}
|g(\zeta)|}{2\pi C^2/2}
\leq C_3,
$$
where $C_3$ is a constant not depending on $z$ and $n$. This implies that
$f(v)=g(0)v^{-1/2}+\alpha(v^{1/2})$ in some neighborhood of $v=0$ cut
along the positive ray.
Let $h$ be so small that $d_1<\frac{C^4}{16}$. Then
$$
\int\limits_{\cal{D}}e^{-\frac{n}{\rho}v}f(v)dv=g(0)\int\limits_{\cal{D}}
e^{-\frac{n}{\rho}v} v^{-1/2} dv+\frac{1}{2\pi i}\int\limits_{\cal{D}}
e^{-\frac{n}{\rho}v}\alpha(v^{1/2})dv
=:g(0)I_1+\frac{1}{2\pi i}I_2 .
$$
Note that
$$
I_2=\int\limits_{\cal{D}}
e^{-\frac{n}{\rho}v}\alpha(v)dv=O\left(\frac{1}{n}\right)
$$
and
$$
I_1=\int\limits_{{\cal{D}}}e^{-\frac{n}{\rho}v}v^{-1/2}dv=
\frac{1}{\left(n/\rho\right)^{1/2}}\int\limits_{\frac{n}{\rho}
\cal{D}}e^{-u}u^{-1/2}du=-\frac{2\Gamma\left(\frac{1}{2}\right)
(1+o(1))}{\left( n/\rho\right)^{1/2}}.
$$
Thus,
$$
\int\limits_{\cal{D}} e^{-\frac{n}{\rho}v}f(v)dv=
\frac{i\sqrt{2\pi}e^{\frac{\rho-1}{2\rho}}(1+o(1)) }{\rho (1-z)
\left(\frac{n}{\rho} \right)^{1/2}}.
$$
$\Box$

\bigskip

\noindent  Lemmas 2 and 3 imply that, as $n\to\infty$,
$$
K_n(z)= \frac{e^{\frac{n}{\rho}}}{R_n^{n+1}}
\frac{i\sqrt{2\pi}e^{\frac{\rho-1}{2\rho}}(1+o(1)) }{\rho (1-z)
\left(n/\rho \right)^{1/2}}, \quad z\in \Omega,
$$
 and hence, by (1.7), (3.1) and Stirling's formula for $\Gamma(x)$, $x>0$,
$$
\frac{\rho (R_n z)^{n+1}}{2\pi i} K_n(z)=
\frac{R_n^nz^n}{\GG}\frac{z}{1-z} (1+o(1)), \quad n\to\infty,
\quad z\in \Omega. \eqno(3.2)
$$
Then Theorem 1 follows from  (2.3)-(2.5), (3.1) and (3.2).

\bigskip

\centerline{{\bf 4. Proof of Theorem 2}}

\bigskip

Denote by
$$
J_1':= \frac{e^{R_n^{\rho} z^{\rho}} \GG }{R_n^n z^{n}} \qquad
\mbox{and} \qquad  J_2':=\frac{\GG }{ R_n^{n+1}z^{n+1}}.
\eqno(4.1)
$$

\bigskip

\noindent{\bf Lemma 4.} {\it Let $z\ne 0$, then for all
sufficiently large $n$,
$$
|J_1'|\geq C_1\left(\frac{n}{\rho}\right)^{1/2},\quad \mbox{if} \quad z\in
\Omega_1\cup\Omega_5, \qquad \mbox{and} \qquad |J_1'|=o(1),\quad
\mbox{if}\quad z\in \Omega_2,
$$
$$
|J_2'|=o(1),\quad \mbox{if}\quad z\in \Omega_3, \qquad
\mbox{and}\qquad |J_2'|\geq (1+(he^{1/\rho})/2)^n, \quad
\mbox{if}\quad z\in \Omega_4,
$$
where $C_1$ is some constant not depending on $n$ and $z$.  }

\bigskip

\noindent {\it Proof}. Using Stirling's formula for $\Gamma(x)$,
$x>0$, and (1.7), as $n\to \infty$,  we have
$$
|J_1'|=\displaystyle{
\frac{e^{\frac{n}{\rho}(1+\frac{\rho-1}{2\rho }+O(1/n^2))r^\rho
\cos(\rho\phi)}\left(\frac{n}{\rho}\right)^{n/\rho+1/2}\sqrt{2\pi}e^{-n/\rho}}{
\left(\frac{n}{\rho}\right)^{n/\rho}e^{\frac{\rho-1}{2\rho}}r^n}(1+o(1))}
$$
$$
=\sqrt{2\pi}(1+o(1))e^{\frac{\rho-1}{2\rho}((1+o(1))r^\rho\cos(\rho\phi)-1)}
\exp\left\{\frac{n}{\rho}\{r^\rho\cos(\rho \phi) -1-\rho\log r\}\right\}
\left(\frac{n}{\rho}\right)^{1/2}
\eqno(4.2)
$$
and
$$
|J_2'|=\frac{\left(\frac{n}{\rho}\right)^{n/\rho+1/2}\sqrt{2\pi}
e^{-n/\rho}}{\left(\frac{n}{\rho}\right)^{n/\rho+1/\rho}
e^{\frac{\rho-1}{2\rho}}r^{n+1}}(1+o(1))=
\frac{\sqrt{2\pi}e^{1/\rho}}{e^{\frac{\rho-1}{2\rho}}}(1+o(1))
\left(\frac{e^{-1/\rho}}{r}\right)^{n+1}\left(\frac{n}{\rho}
\right)^{1/2-1/\rho}.
\eqno(4.3)
$$
Equations (4.2) and (4.3) make Lemma 4 obvious.   $\Box$

\bigskip

According to Theorem 1 and (4.1), we have
$$
\frac{I_n(R_n z;\lambda)\GG}{R_n z^n}=-\lambda \rho (1+o(1))J_1'-
(1+o(1))\frac{z}{1-z}, \qquad z\in \Omega_5\cup\left\{
\Omega_2\cap\{z: |z|\geq 1\}\right\};
$$
$$
\frac{I_n(R_n z;\lambda)\GG}{R_n z^n}=(1-\lambda) \rho (1+o(1))J_1'-
(1+o(1))\frac{z}{1-z}, \qquad z\in \Omega_1\cup\left\{
\Omega_2\cap\{z: |z|\leq 1\}\right\};
$$
$$
\frac{I_n(R_n z;\lambda)\GG}{R_n^n z^n}=\frac{(\lambda-1)}{\G}J_2'-
(1+o(1))\frac{z}{1-z}, \quad z\in \Omega_3\cup\Omega_4.
$$
The last three expressions and Lemma 4 completes the proof of Theorem 2.

\bigskip

\centerline{{\bf 5. Proof of Theorem 3}}

\bigskip

We can rewrite $I_n(R_n z;\lambda)$ as
$$
I_n(R_n z;\lambda)=\ss-\lambda E_{\frac{1}{\rho}}(R_n z).
\eqno(5.1)
$$
Theorem 1 of \cite{[4]} implies that
$$
\frac{ s_n \left( R_n \left(1+\left(\frac{2}{\rho n}
\right)^{1/2}\zeta\right);E_{\frac{1}{\rho}} \right)}{ \left(1+
\left(\frac{2}{\rho n} \right)^{1/2}\zeta \right)^n
E_{\frac{1}{\rho}} \left( R_n\right)}\to
\frac{1}{2}e^{\zeta^2}\mbox{erfc}(\zeta), \eqno(5.2)
$$
as $n\to\infty$, uniformly on every compact set of $\zeta$ -- plane.

By (2.2), we have
$$
\frac{ E_{\frac{1}{\rho}} \left( R_n \left(1+\left(\frac{2}{\rho
n} \right)^{1/2}\zeta\right);E_{\frac{1}{\rho}} \right)}{ \left(1+
\left(\frac{2}{\rho n} \right)^{1/2}\zeta \right)^n
E_{\frac{1}{\rho}} \left( R_n\right)}= \exp
\left\{\left(\frac{n}{\rho}+\frac{\rho-1}{2\rho}\right)
\left(1+\left(\frac{2}{\rho n} \right)^{1/2}\zeta \right)^\rho
\right.
$$
$$
\left.
-n\log \left(1+\left(\frac{2}{\rho n}\right)
^{1/2}\zeta\right)-\left(\frac{n}{\rho}+\frac{\rho-1}{2\rho}\right)+o(1)
\right\} =e^{\zeta^2}+o(1),\quad n\to\infty. \eqno(5.3)
$$
Theorem 3 follows from (5.1) -- (5.3).

\bigskip

\centerline {{\bf 6. Proof of Theorem 4}}

\bigskip

\noindent I. By Theorem 1, as $n\to \infty$, we have
$$
I_n\left(R_n \xi \left(1+\frac{\log n}{2(1-\xi^\rho)n}-
\frac{\zeta- i\tau_n}{(1-\xi^\rho)n}\right);\lambda\right)
\displaystyle{ \frac{\GG}{R_n^n \xi^n\left(1+\frac{\log
n}{2(1-\xi^\rho)n}- \frac{\zeta- i\tau_n}{(1-\xi^\rho)n}\right)^n}
}
$$
$$
=\displaystyle{\frac{\mu\rho e^{R_n^\rho
\xi^\rho\left(1+\frac{\log
n}{2(1-\xi^\rho)n}-\frac{\zeta-i\tau_n}{(1-\xi^\rho)n}\right)^\rho}\GG}{R_n^n
\xi^n\left(1+\frac{\log
n}{2(1-\xi^\rho)n}-\frac{\zeta-i\tau_n}{(1-\xi^\rho)n}\right)^n}}-\frac{\xi}{1-\xi}(1+o(1))
$$
$$
=:\mu\rho A_n-\frac{\xi}{1-\xi}(1+o(1)),
$$
where $\mu=1-\lambda$, if $|\xi|<1$, and $\mu=-\lambda$, if
$|\xi|>1$. Then, the first assertion of the Theorem follows from
the following equality:
$$
A_n=(2\pi \rho^{-1})^{1/2} \exp \left\{
\frac{n}{\rho}(|\xi|^\rho\cos \rho \phi -1 -\rho \log
|\xi|)+i\frac{n}{\rho}(|\xi|^\rho \sin(\rho \phi)-\rho\phi)
-i\tau_n +\zeta\right.
$$
$$
\left. +\frac{\rho+1}{2\rho}(\xi^\rho -1)+O\left(\frac{\log^2
n}{n}\right)\right\}=\displaystyle{(2\pi
\rho^{-1})^{1/2}e^{\frac{\rho+1}{2\rho}(\xi^\rho-1)} e^{\zeta}
+o(1)}, \quad n\to\infty,
$$
which one can easily obtain using Stirling's formula for
$\Gamma(1+x)$, $x>0$, and (1.7).

\bigskip

\noindent II. By Theorem 1, as $n\to\infty$, we have
$$
I_n\left(R_n \xi \left(1+\left(\frac{1}{2}-\frac{1}{\rho}\right)
\frac{\log n}{n}- \frac{\zeta-
i\tau^{'}_n}{n+1}\right);\lambda\right) \displaystyle{
\frac{\GG}{R_n^n \xi^n\left(1+(\frac{1}{2}-\frac{1}{\rho})
\frac{\log n}{n}- \frac{\zeta- i\tau^{'}_n}{n+1}\right)^n}}
$$
$$
=\displaystyle{\frac{(\lambda -1) \GG}{\G R_n^{n+1} \xi^{n+1}
\left(1+\left(\frac{1}{2}-\frac{1}{\rho}\right)\frac{\log n}{n}
-\frac{\zeta -i\tau'_n
}{n+1}\right)^{n+1}}}-\frac{\xi}{1-\xi}(1+o(1))
$$
$$
=:
\frac{\lambda-1}{\G} B_n-\frac{\xi}{1-\xi}(1+o(1)).
$$
The second assertion of the Theorem follows then from the
following equality:
$$
B_n=\left(2\pi  e^{\frac{1-\rho}{\rho}}\right)^{1/2}
\rho^{\frac{1}{\rho}-\frac{1}{2}}\exp\left\{-\zeta+i\tau'_n-i(n+1)\phi
+O\left(\frac{1}{n}\right)\right\}
$$
$$
=\left(2\pi
e^{\frac{1-\rho}{\rho}}\right)^{1/2}
\rho^{\frac{1}{\rho}-\frac{1}{2}}e^{-\zeta}+o(1), \quad
n\to\infty,
$$
which is easy to check using Stirling's formula for $\Gamma(1+x)$,
$x>0$, and (1.7). $\Box$

\bigskip

\noindent{\bf Acknowledgments}: The author is grateful to Prof.
C.Y. Y{\i}ld{\i}r{\i}m for constant attention to this work and to
Prof. I.V. Ostrovskii for useful discussions.

\end{large}

\bigskip

\noindent
Department of Mathematics,\\
Bilkent University,\\
06533 Bilkent, Ankara, Turkey,\\
E-mail: natalya@fen.bilkent.edu.tr

\end{document}